\documentclass[12pt,a4paper]{article}
\usepackage{a4wide,amsmath,amsthm,amssymb,latexsym}
\usepackage{enumerate}

\newtheorem{theorem}{Theorem}
\newtheorem{lemma}{Lemma}
\newtheorem{corollary}{Corollary}
\newtheorem{proposition}{Proposition}

\newtheorem{remark}{Remark}

\newcommand{\Z}{{\mathbb Z}}
\newcommand{\N}{{\mathbb N}}

\newcommand{\R}{{\mathbb R}}

\newcommand{\IS}{{\mathcal I}{\mathcal S}}

\newcommand{\BB}{{\mathfrak B}}

\renewcommand{\phi}{\varphi}

\newcommand{\deff}{\operatorname{def}}
\newcommand{\ran}{\operatorname{im}}
\newcommand{\rank}{\operatorname{rank}}
\newcommand{\str}{\operatorname{st.rank}}
\newcommand{\dom}{\operatorname{dom}}

\renewcommand{\Pr}{\operatorname{Pr}}

\begin{document}

\title{Combinatorics and distributions of partial injections}
\author{Olexandr Ganyushkin and Volodymyr Mazorchuk}
\date{}
\maketitle

\begin{abstract}
We obtain several combinatorial results about chains, cycles and
orbits of the elements of the symmetric inverse semigroup $\IS_n$ 
and the set $T_n$ of nilpotent elements in $\IS_n$. We also get  
some estimates for  the growth of $|\IS_n|$ and $|T_n|$, and study 
random products of elements from $\IS_n$.
\end{abstract}

\section{Introduction}\label{s1}

Roughly speaking, there are three semigroups, which play a principal
role in the theory of transformation semigroups. The first one is the
{\em full transformation semigroup} $\mathcal{T}_M$ of all
transformations of a set, $M$, the second one is the {\em full partial
transformation semigroup} $\mathcal{PT}_M$ of all partial
transformations of $M$, and the third one is the {\em symmetric
inverse semigroup} $\IS_M$ of all partial injective transformations of 
$M$. The role of the last semigroup is especially important in the
theory of inverse semigroups, where this role is  analogous to that of
the symmetric group $S_n$ in the group theory.  

In the present paper we consider only finite sets. Hence we choose 
$M$ to be the set $N=\{1,2,\dots,n\}$. We denote the corresponding
semigroups by $\mathcal{T}_n$, $\mathcal{PT}_n$ and $\IS_n$
respectively. 

Combinatorial properties of $\mathcal{T}_n$, $\mathcal{PT}_n$ and some 
related transformation semigroups (for example the semigroup 
$\mathcal{O}_n=\{\alpha\in\mathcal{T}_n:x\leq
y\Rightarrow\alpha(x)\leq \alpha(y)\}$ of all transformations
preserving the natural order) were studied in a number of papers by
several authors, see for example \cite{Hi1,Hi2,Ho,GH1,Ka} and
references therein. In particular, in the monograph \cite{Hi3} many 
combinatorial properties of $\mathcal{T}_n$ are collected in a 
separate big chapter. At the same time the situation with the
semigroup $\IS_n$ is completely different. Only few papers, in which 
nilpotent elements and nilpotent subsemigroups of $\IS_n$ are studied, 
deal partially with some combinatorial questions, see 
\cite{GH2,GK1,GK2,GP,GM}. A survey on these combinatorial results and 
some new combinatorial results on $\IS_n$ can be found in 
\cite{GM2}. Monographs on semigroups, even those, dedicated
completely to inverse semigroups, for example \cite{La,Pe}, do not
go much further than giving a formula for the cardinality of $\IS_n$. 
No combinatorial results can be found even in the monograph \cite{Li}, 
which is dedicated completely to $\IS_n$.  

In the present paper we study combinatorial properties of the elements
of $\IS_n$ in general. The action of the element $\alpha\in\IS_n$ on
$N$ is described by the graph of the action, which leads to the
standard combinatorial data, including such notions as {\em
cyclic} and {\em chain} components of the graph and {\em orbits} of
the elements from $N$. In Sections~\ref{s2} and \ref{s3} we obtain
several combinatorial formulae relating the ingredients of these data
with each other, with the cardinality of the semigroup $\IS_n$
itself, and with the cardinality of the set $T_n$ of all nilpotent
elements from $\IS_n$. 

In Section~\ref{s4} we concentrate on the study of the set $T_n$ and
discover possibly the most surprising result of the paper, namely a
strange duality between $T_n$ and $\IS_n$. This duality is
incarnated into a number of statements, each consisting of a pair of
equalities, dual to each other in the sense, that one of the
equalities is obtained from the other one by substituting the 
combinatorial data, related  to $\IS_n$, with the corresponding 
combinatorial data, related to $T_n$, and vice versa.

In Section~\ref{s5} we study the asymptotics of both $|\IS_n|$ and 
$|T_n|$. We obtain that the growth of both $|\IS_n|$ and $|T_n|$
can be (very) roughly described by $(n+2)!$, in particular, that it is
roughly the same. At the same time, it is also shown that the limit
value of the ratio $|T_n|/|\IS_n|$ is $0$.

Finally, in Section~\ref{s6} we study random products of $k$ elements
from $\IS_n$ under the assumption of the uniform distribution of original
probabilities. We give both, a precise formula and some estimates, for 
the probability of such product to equal some fixed element from $\IS_n$,
and show that for all $k$ big enough almost all products of $k$
elements from $\IS_n$ are zero. The distribution of probabilities we
calculate is controlled by a square upper triangular matrix with 
non-negative integer entries. We show that the eigenvectors of this
matrix can be computed purely combinatorially, in terms of the
combinatorial data of $\IS_n$, and derive that the corresponding
transformation matrix transforms the vector $(1,1,\dots,1)$ into the
vector $(|\IS_n|,|\IS_{n-1}|,\dots,|\IS_0|)$.

\section{Preliminary combinatorics}\label{s2}

Throughout the paper for two sets, $X$ and $Y$, by $X\subset Y$ we 
mean that $x\in X$ implies $x\in Y$ for every element $x$ (in particular,
$X=Y$ implies $X\subset Y$).

From the definition of $\IS_n$ it follows immediately that every
element $a\in\IS_n$ is uniquely determined by its domain $\dom(a)$,
its range $\ran(a)$ and a bijection from $\dom(a)$ to $\ran(a)$. Hence
\begin{displaymath}
|\IS_n|=\sum_{k=0}^{n}\binom{n}{k}^2 k!.
\end{displaymath}
The number $\rank(a)=|\dom(a)|=|\ran(a)|$ is called the {\em rank} of
$a$ and the number $\deff(a)=n-\rank(a)$ is called the {\em defect} of
$a$. 

For elements from $\IS_n$ one can use their regular table
presentation
\begin{displaymath}
a=\left(\begin{array}{cccc} i_1 & i_2 & \dots & i_k \\
j_1 & j_2 & \dots & j_k \end{array}\right),
\end{displaymath}
where $\dom(a)=\{i_1,\dots,i_k\}$ and
$\ran(a)=\{j_1,\dots,j_k\}$. However, sometimes it is more convenient
to use the so-called {\em chain} (or {\em chart}) decomposition of
$a$, which is analogous to the cyclic decomposition for usual
permutations. We refer the reader to \cite{Li} for rigorous
definitions, however, this decomposition is very easy to explain on the
following example. The element
\begin{displaymath}
a=\left(\begin{array}{ccccccc}
1 & 2 & 3 & 4 & 5 & 7 & 9 \\7 & 4 & 5 & 1 & 10 & 2 & 6 
\end{array}
\right)\in\IS_{10}
\end{displaymath}
has the following graph of the action on $\{1,2,\dots,10\}$:
\begin{displaymath}
\begin{array}{ccc}
1 & \to & 7 \\  \uparrow & & \downarrow \\ 4 & \leftarrow & 2
\end{array}
\quad\quad
3\to 5\to 10\quad\quad 9\to 6 \quad\quad 8,
\end{displaymath}
and hence it is convenient to write it as $a=(1,7,2,4)[3,5,10][9,6][8]$. We
call $(1,7,2,4)$ a {\em cycle} and $[3,5,10]$ (as well as $[9,6]$ and $[8]$)
a {\em chain} of the element $a$. We remark that chains of length $1$
correspond to those elements $x\in N$, which do not belong to
$\dom(a)\cup\ran(a)$. It is obvious that $\deff(a)$ equals the number
of chains in the chain decomposition of $a$. For $a\in\IS_n$ and $i\in N$ let
$c_i$ and $d_i$ denote respectively the number of cycles and the number of chains 
of length $i$ in the chain decomposition of $a$. The vector 
$(c_1,\dots,c_n,d_1,\dots,d_n)$ is called the {\em chain type} of $a$, see 
\cite{GM,Li}.

\begin{proposition}\label{p1}(\cite[Lemma~V.1.9]{Ho1})
The set $E(\IS_n)$ of idempotents in $\IS_n$ is a semigroups,
isomorphic to the semigroup $\BB_n=\{A\,:\, A\subset N\}$ with the
intersection of sets as the corresponding binary operation. In 
particular, $|E(\IS_n)|=2^n$. 
\end{proposition} 

The semigroup $\IS_n$ contains the zero element $0$, which is the
unique transformation  such that $\dom(0)=\varnothing$. Recall that if
$S$ is a semigroup with zero $0$, then the element $a\in S$ is called
{\em nilpotent} provided that $a^k=0$ for some $k>0$. We will denote by
$T_n$ the set of all nilpotent elements in $\IS_n$ and remark that
$T_n$ is not a subsemigroup of $\IS_n$ (the product of two nilpotent
elements is not nilpotent in general). 

\begin{proposition}\label{p2}(\cite{GM2})
The element $a\in\IS_n$ is nilpotent if and only if the chain
decomposition of $a$ contains only chains. The number of nilpotent
elements in $\IS_n$ with the given defect $k$ equals the signless Lah
number $L'(n,k)=\frac{n!}{k!}\binom{n-1}{k-1}$.
\end{proposition} 

By the {\em permutational part} of the element 
$(a_1,\dots,a_p)\dots(b_1,\dots,b_q)[c_1,\dots,c_s][d_1,\dots,d_t]$ we
will mean the element $(a_1,\dots,a_p)\dots(b_1,\dots,b_q)
[c_1][c_2]\dots[d_t]$. The rank of the permutational part of
$\alpha\in\IS_n$ is called the {\em stable rank} of $\alpha$ and is
denoted by $\str(\alpha)$. This notion is analogous to the
corresponding notion for $\mathcal{T}_n$, see \cite{Hi3}. It is
obvious that $\str(\alpha)=\str(\alpha^i)$ for all $i\in\N$.

Taking the inverse element defines an anti-involution, 
$\alpha\mapsto\alpha^{-1}$, on $\IS_n$. The action of this 
anti-involution on $\alpha$ can be described as follows: one takes the
graph of the action of $\alpha$ and reverses all arrows in it. It
follows that this map does not change the chain type of $\alpha$, in 
particular, nilpotent elements are sent to nilpotent elements. Since 
this map switches $\ran(\alpha)$ and $\dom(\alpha)$, it allows one to 
transfer all statements about the ranges of the elements (in 
particular, of nilpotents) to the dual statements about the domains, 
and vice versa. 

Studying probabilistic characteristics of various parameters of
elements in $\IS_n$ it is natural to assume that the original
distribution of probabilities of the elements in $\IS_n$ is
uniform. An unexpected difficulty in this case is the fact that 
for two fixed $x,y\in N$ the
random events ``$x\in\dom(\alpha)$'' and ``$y\in\dom(\alpha)$'' are
not independent in general. For example, in $\IS_3$ we have 
\begin{gather*}
\Pr\big(1\in \dom(\alpha)\big)=\Pr\big(2\in \dom(\alpha)\big)=21/34,
\quad\text{but}\\
\Pr\big((1\in \dom(\alpha))\text{ and }(2\in \dom(\alpha))\big)=
6/17\neq (21/34)^2.
\end{gather*}
Furthermore, the random events ``$x\in\dom(\alpha)$'' and 
``$y\in\dom(\alpha)$'' are not independent if we consider them for 
$T_n$ instead of $\IS_n$ either. 

For $k=0,\dots,n$ denote
\begin{gather*}
R_{n,k}=|\{\alpha\in\IS_n:\rank(\alpha)=k\}|,\\ 
D_{n,k}=|\{\alpha\in\IS_n:\deff(\alpha)=k\}|,\\ 
St_{n,k}=|\{\alpha\in\IS_n:\str(\alpha)=k\}|.
\end{gather*}
Then we have
\begin{displaymath}
R_{n,k}=\binom{n}{k}^2\cdot k!\quad\quad\text{ and }\quad\quad
|\IS_n|=\sum_{k=0}^n R_{n,k}.
\end{displaymath}
As $\rank(\alpha)+\deff(\alpha)=n$, we have 
\begin{displaymath}
D_{n,k}=R_{n,n-k}=
\binom{n}{k}^2\cdot (n-k)!\quad\quad\text{ and }\quad\quad
|\IS_n|=\sum_{k=0}^n D_{n,k}.
\end{displaymath}
From Proposition~\ref{p2} we have $\displaystyle|T_n|=\sum_{k=1}^n
L'(n,k)$. Now from 
\begin{displaymath}
L'(n,k)=\binom{n-1}{k-1}\frac{n!}{k!}=\frac{k}{n}\binom{n}{k}^2(n-k)!=
\frac{k}{n}D_{n,k}
\end{displaymath}
it follows that 
\begin{equation}\label{eqeq1}
|T_n|=\sum_{k=1}^n\frac{k}{n}D_{n,k}=\sum_{k=1}^n\frac{n-k}{n}R_{n,k}.
\end{equation}

\begin{remark}\label{r2-3}{\rm
There is a purely combinatorial way to show that the sets
\begin{gather*}
M_1=\{(\alpha,x):\alpha\in T_n,\deff(\alpha)=k,x\in N\}\quad \text{
and }\\
M_2=\{(\beta,l):\beta\in\IS_n,\deff(\beta)=k,l\text{ is a chain 
of }\beta\} 
\end{gather*}
have the same cardinality, which 
implies $nL'(n,k)=kD_{n,k}$. Indeed, for $(\alpha,x)\in M_1$ we define 
$f((\alpha,x))=(\beta,l)\in M_2$ in the following way: let
$N_{\alpha}=\{y\in N:\alpha^i(x)=y\text{ for some }i\in \N\}=
\{t_1,t_2,\dots,t_s\}$, $t_1<t_2,\dots<t_s$, then
$\dom(\beta)=(\dom(\alpha)\cup N_{\alpha})\setminus\{x\}$,
$\beta(y)=\alpha(y)$ for all $y\in\dom(\alpha)\setminus
(\{x\}\cup N_{\alpha})$, $\beta(\alpha^i(x))=t_i$ for all 
$i=1,\dots,s$; and $l$ is the chain of $\beta$, containing $x$. One 
easily checks that $(\beta,l)\in M_2$ and that $f$ is a bijection. 
}
\end{remark}

For $x\in \R$ and $k\in\{0,1,\dots\}$ we denote by $[x]_k$ that {\em
$k$-th decreasing factorial} $[x]_k=x(x-1)\dots(x-k+1)$. 

\begin{proposition}\label{p3}
$\displaystyle St_{n,k}=[n]_k\sum_{i=1}^{n-k}L'(n-k,i)$.
\end{proposition} 

\begin{proof}
We partition $\IS_n$ into classes with respect to the domain $A\subset
N$ of the permutational part of the element $\alpha\in\IS_n$. The
element $\alpha$ acts as a permutation on $A$ and as a nilpotent on
$N\setminus A$. Choosing $A$ such that $|A|=k$, a permutation on $A$, 
and a nilpotent on $N\setminus A$ in all possible ways, and taking
Proposition~\ref{p2} into account, we get
\begin{displaymath}
St_{n,k}=\binom{n}{k}\cdot k!\cdot
\sum_{i=1}^{n-k}L'(n-k,i)=[n]_k\sum_{i=1}^{n-k}L'(n-k,i).
\end{displaymath}
\end{proof}

Denote by $C_{n,k}$ the number of cycles of length $k$ and by
$L_{n,k}$ the number of all chains of length $k$ in all elements 
in $\IS_n$.

\begin{proposition}\label{p4}
$L_{n,k}=[n]_k\cdot|\IS_{n-k}|$ and 
$C_{n,k}=\frac{1}{k}[n]_k\cdot|\IS_{n-k}|$.
\end{proposition} 

\begin{proof}
The number of those elements in $\IS_n$, whose chain decomposition
contains a fixed cycle (chain) of length $k$, equals $|\IS_{n-k}|$. On
the other hand, given $k$ elements from $N$, we can form $k!$ different
chains and $(k-1)!$ different cycles of length $k$. Now the remark
that $\binom{n}{k}\cdot k!=[n]_k$ completes the proof.
\end{proof}

Invertible elements in $\IS_n$ are exactly permutations, that is
elements of the symmetric group $S_n$. Hence
$b_n=\frac{|\IS_n|}{|S_n|}=\frac{|\IS_n|}{n!}$ characterizes  (in some
sense) the non-invertability of the elements of $\IS_n$, or how far
$\IS_n$ is from being a group.

\begin{corollary}\label{c1}
The average number $c_n$ of components in the chain decomposition of
the element $\alpha\in\IS_n$ equals
\begin{displaymath}
c_n=b_n^{-1}\sum_{k=1}^n\left(1+\frac{1}{k}\right)b_{n-k}.
\end{displaymath}
\end{corollary} 

\begin{proof}
From Proposition~\ref{p4} it follows that
\begin{displaymath}
c_n=\frac{1}{|\IS_n|}
\sum_{k=1}^n[n]_k\left(1+\frac{1}{k}\right)|\IS_{n-k}|
=\frac{n!}{|\IS_n|}\sum_{k=1}^n\left(1+\frac{1}{k}\right)
\frac{|\IS_{n-k}|}{(n-k)!}.
\end{displaymath}
\end{proof}

\begin{remark}\label{r6}{\rm
One can compare the last result with $S_n$ and $\mathcal{T}_n$:
the average number of components (that is cycles) in the cyclic
decomposition of a permutation $\pi\in S_n$ equals $1+\frac{1}{2}+
\frac{1}{3}+\cdots +\frac{1}{n}$, and the average number of components
for an element $f\in \mathcal{T}_n$ equals $\displaystyle
\sum_{k=1}^n\frac{n!}{k(n-k)!n^k}$, see \cite[Lemma~6.1.12]{Hi3}
or \cite{Kr}.
}
\end{remark}
 
\section{Chains and Orbits}\label{s3}

Let $L_n$ denote the total number of chains in the chain
decompositions of all elements in $\IS_n$.

\begin{proposition}\label{p5}
$\displaystyle L_n=\sum_{k=0}^n (n-k)R_{n,k}$.
\end{proposition}

\begin{proof}
Each element of rank $k$ has defect $n-k$ and thus contains $n-k$ chains.
\end{proof}

Comparing the last formula with Proposition~\ref{p4} we get

\begin{corollary}\label{c2}
\begin{displaymath}
\sum_{k=0}^{n-1} (n-k)\binom{n}{k}^2k!=\sum_{k=1}^{n}[n]_k|\IS_{n-k}|.
\end{displaymath}
\end{corollary}

One more recursive relation for the cardinalities of $\IS_n$ is given
by

\begin{proposition}\label{p6}
\begin{displaymath}
\frac{1}{n}
\sum_{k=1}^{n}\big(k\cdot R_{n,k}+[n]_k|\IS_{n-k}|\big)=|\IS_n|.
\end{displaymath}
\end{proposition}

\begin{proof}
We have $\rank(\alpha)+\deff(\alpha)=n$ for every
$\alpha\in\IS_n$. Hence the sum of the average rank and the average 
defect of all elements in $\IS_n$ must be 
equal to $n$ as well. Therefore
\begin{displaymath}
\frac{1}{|\IS_n|}\sum_{k=1}^{n}k\cdot R_{n,k}+
\frac{1}{|\IS_n|}\sum_{k=1}^{n}
[n]_k|\IS_{n-k}|=n,
\end{displaymath}
which completes the proof.
\end{proof}

\begin{theorem}\label{t2}
\begin{displaymath}
L_n=\sum_{\alpha\in \IS_n}\str(\alpha).
\end{displaymath}
\end{theorem}

\begin{proof}
Consider the sets
\begin{gather*}
A=\{(\alpha,c,x):\alpha\in\IS_n,c \text{ is a cycle of }\alpha,
x\text{ is a point of }c\},\\
B=\{(\beta,l):\beta\in\IS_n,l \text{ is a chain from the chain
decomposition of }\beta\}.
\end{gather*}
The statement of the theorem is equivalent to the equality $|A|=|B|$. 

Consider the map $f:A\to B$, which is defined as follows:
$f((\alpha,(x,a,\dots,b),x))=(\beta,[x,a,\dots,b])$, where $\beta$
is obtained from $\alpha$ substituting  the cycle $(x,a,\dots,b)$ with
the chain $[x,a,\dots,b]$. Consider also the map $g:B\to A$, which is 
defined as follows: $g((\beta,[x,a,\dots,b]))=(\alpha,(x,a,\dots,b),x)$, 
where $\alpha$ is obtained from $\beta$ substituting the chain
$[x,a,\dots,b]$ with the cycle $(x,a,\dots,b)$. Obviously  $f$ and 
$g$ are inverse to each other and thus $|A|=|B|$.
\end{proof}

\begin{remark}\label{r3}{\rm
It is obvious that $\displaystyle\sum_{\alpha\in \IS_n}\str(\alpha)$
is equal to the total sum of lengths of all cycles of all elements in
$\IS_n$. 
}
\end{remark}

Let $P_n$ denote the total number of fixed points for all elements in
$\IS_n$. From Burnside's lemma it follows that the average number of
fixed points for permutations in $S_n$ equals $1$. An analogue of this
statement for $\IS_n$ is the following 

\begin{theorem}\label{t1}
\begin{equation}\label{eq1}
P_n+\frac{1}{n}L_n=|\IS_n|.
\end{equation}
\end{theorem}

\begin{proof}
Consider the following sets:
\begin{gather*}
A=\{(\alpha,x):\alpha\in\IS_n,x\in N\},\\
B=\{(\beta,l):\beta\in\IS_n, l \text{ is a chain for the chain
decomposition of }\beta\},\\
C=\{(\gamma,y,z):\gamma\in\IS_n, y\text{ is a fixed point of
}\gamma,z\in N\}.
\end{gather*}
The equality \eqref{eq1} is equivalent to the equality
$|A|=|B|+|C|$. To prove the latter we decompose $A$ into a disjoint
union $A=A_1\cup A_2$, where
\begin{gather*}
A_1=\{(\alpha,x)\in A:x\text{ belongs to some chain of }\alpha\},\\
A_2=\{(\alpha,x)\in A:x\text{ belongs to some cycle of }\alpha\}.
\end{gather*}
Consider the transformation, which maps the cycle $(x,a,\dots,b)$ with
a base point $x$ to the chain $[x,a,\dots,b]$. Obviously, this
transformation induces a bijection $A_2\to B$. Hence $|A_2|=|B|$.

To prove $|A_1|=|C|$ we construct mutually inverse bijections
$f:A_1\to C$ and $g:C\to A_1$. Consider any element $(\alpha,x)\in
A_1$. If $x$ is the source of some chain $[x,a,\dots,b]$ of length 
at least $2$ from the chain decomposition of $\alpha$, we define
$f((\alpha,x))=(\gamma,x,a)$, where $\gamma$ is obtained from 
$\alpha$ substituting the chain $[x,a,\dots,b]$ with the cycle
$(x)$ and the cycle 
$(a,\dots,b)$. If $x$ is the only point of the chain $[x]$, we 
define $f((\alpha,x))=(\gamma,x,x)$, where $\gamma$ is obtained from 
$\alpha$ substituting the chain $[x]$ with the cycle $(x)$. 
Finally, if $x$ is contained in some chain $[a,\dots,b,x,c,\dots,d]$ 
and is different from the source of this chain, we define 
$f((\alpha,x))=(\gamma,x,b)$, where $\gamma$ is obtained from 
$\alpha$ substituting the chain $[a,\dots,b,x,c,\dots,d]$ with 
the cycle $(x)$ and the chain $[a,\dots,b,c,\dots,d]$. 

Let now $(\gamma,y,z)\in C$. If $y=z$, we define
$g((\gamma,y,z))=(\alpha,z)$, where $\alpha$ is obtained from $\gamma$
substituting the cycle $(y)$ with the chain $[y]$. If $z$ is  a
point of some chain $[a_1,\dots,a_s,z,b_1,\dots,b_t]$ in the chain
decomposition of $\gamma$, we set $g((\gamma,y,z))=(\alpha,z)$, where
$\alpha$ is obtained from $\gamma$ substituting the cycle $(y)$ and
the chain $[a_1,\dots,a_s,z,b_1,\dots,b_t]$ with the chain 
$[a_1,\dots,a_s,z,y,b_1,\dots,b_t]$. Finally, if $z$ is  a
point of some cycle $(a_1,\dots,a_s,z)$ of $\gamma$, we set 
$g((\gamma,y,z))=(\alpha,z)$, where $\alpha$ is obtained from 
$\gamma$ substituting the cycles $(y)$ and $(z,a_1,\dots,a_s)$ with 
the chain $[y,z,a_1,\dots,a_s]$.

Obviously, $f$ and $g$ are inverse to each other implying $|A_1|=|C|$,
and the theorem follows.
\end{proof}

If $x\in\dom(\alpha)$, the set
$\{x,\alpha(x),\alpha^2(x),\dots\}$ is called the {\em orbit} of
$x$ under $\alpha$ and the cardinality of this set is called the {\em
length} of the orbit. If $x\not\in \dom(\alpha)$, we say that the
orbit is empty and consequently the length of the orbit is $0$. Since
for every transposition $(x,y)\in S_n$ the conjugation $\alpha\mapsto
(x,y)\alpha(x,y)$ maps orbits of $x$ to orbits of $y$ and vice versa,
it is enough to study the orbits of the element $1$.

It is easy to calculate that the average length of the orbit of $1$
under the action of the symmetric group $S_n$ equals $(n+1)/2$, and 
the number of orbits of $1$ of length $i$ does not depend on $i$ and 
equals $(n-1)!$. The corresponding situation in the semigroup
$\mathcal{T}_n$ is much more interesting. For example, it is shown in
\cite{Ha} that the random function $X_n(\alpha)$, whose value is the
cardinality of the permutational part of $\alpha\in \mathcal{T}_n$,
and the random function $Y_n(\alpha)$, whose value is the length of
the orbit of $1$ for  $\alpha\in \mathcal{T}_n$, have the same
distribution. Later on an elementary proof of this fact was found in
\cite{BH} (see also the historical review of this fact in \cite{Hi2}).
For $\IS_n$ the corresponding statement does not hold, however, one
has the following

\begin{theorem}\label{t3}
The sum of lengths of the orbits of $1$ over all elements $\alpha\in
\IS_n$ equals the total number of chains in all elements in $\IS_n$.
\end{theorem}

\begin{proof}
Let
\begin{gather*}
A=\{(\alpha,x):\alpha\in\IS_n,x \text{ is a member of the orbit of $1$
for }\alpha\},\\
B=\{(\beta,l):\beta\in\IS_n,l \text{ is a chain from the chain
decomposition of }\beta\}.
\end{gather*}
The statement of the theorem is equivalent to the equality
$|A|=|B|$. To prove the latter let us construct mutually inverse
bijections $f:A\to B$ and $g:B\to A$.

Let $(\alpha,x)\in A$. If $x$ is a point of the cycle
$(x,\dots,1,\dots,y)$, we define
$f((\alpha,x))=(\beta,[x,\dots,1,\dots,y])$, where $\beta$
is obtained from $\alpha$ substituting the cycle 
$(x,\dots,1,\dots,y)$ with the chain $[x,\dots,1,\dots,y]$. If $x$ is
a point of the chain $[a,\dots,1,\dots,b,x,\dots,c]$ and $x\neq 1$, we 
define $f((\alpha,x))=(\beta,[x,\dots,c])$, where $\beta$
is obtained from $\alpha$ substituting the chain 
$[a,\dots,1,\dots,b,x,\dots,c]$ with two chains, $[a,\dots,1,\dots,b]$ 
and $[x,\dots,c]$. Finally, if $x=1$ and it is a point of the chain 
$[a,\dots,1,b,\dots,c]$, we define
$f((\alpha,x))=(\beta,[b,\dots,c])$, where $\beta$ is obtained 
from $\alpha$ substituting the chain 
$[a,\dots,1,b,\dots,c]$ with the cycle $(a,\dots,1)$ and the chain 
$[b,\dots,c]$.

Let $(\beta,l)\in B$. If $l$ contains $1$ and has the form
$l=[x,\dots,1,\dots,y]$, we set $g((\beta,l))=(\alpha,x)$, where  
$\alpha$ is obtained from $\beta$ substituting the chain 
$[x,\dots,1,\dots,y]$ with the cycle  $(x,\dots,1,\dots,y)$. If 
$l=[x,\dots,c]$ does not contain $1$ and $1$ belongs to another 
chain, $[a,\dots,1,\dots,b]$ say, we set $g((\beta,l))=(\alpha,x)$, 
where $\alpha$ is obtained from $\beta$ substituting 
the chains $[a,\dots,1,\dots,b]$ and $[x,\dots,c]$ with 
the chain $[a,\dots,1,\dots,b,x,\dots,c]$. 
If $l=[x,\dots,c]$ does not contain $1$ and $1$ belongs to a cycle, 
$(a,\dots,1)$ say, we set $g((\beta,l))=(\alpha,1)$, where  
$\alpha$ is obtained from $\beta$ substituting the chain 
$[x,\dots,c]$ and the cycle $(a,\dots,1)$ with the chain 
$[a,\dots,1,x,\dots,c]$. It is obvious that under the definition of
$g$ the point $x$ in the pair $(\alpha,x)$ always belongs to the orbit 
of $1$ under the action of $\alpha$. 

It is easy to check that $f$ and $g$ are mutually inverse bijections,
which completes the proof.
\end{proof}

\begin{theorem}\label{t4}
Let $l_{n,k}$ denote the total number of orbits of $1$, having length
$k$, in all elements from $\IS_n$. Then 
\begin{enumerate}[(i)]
\item\label{ff1} $l_{n,0}=|T_n|$,
\item\label{ff2} $l_{n,1}=|\IS_{n-1}|$,
\item\label{ff3} $l_{n,k}=[n-1]_{k-1}(L_{n-k}+2|\IS_{n-k}|)$ for 
$1<k\leq n$.
\end{enumerate}
\end{theorem}

\begin{proof}
\eqref{ff1}. According to the definition, $l_{n,0}$ is the cardinality 
of the set 
\begin{displaymath}
E(n,0)=\{\alpha\in\IS_n:1\not\in\dom(\alpha)\}. 
\end{displaymath}
Consider the following decomposition of $E(n,0)$ into a disjoint union
of subsets:
\begin{displaymath}
E(n,0)=\bigcup_{A\subset\{2,3,\dots,n\}} E_A,
\end{displaymath}
where $\displaystyle E_A=\{\alpha\in
E(n,0):A=\cap_{k>0}\dom(\alpha^k)\}$. In other words, $E_A$ contains
all those elements from $E(n,0)$, for which $A$ is the domain of the
permutational part. 

Consider also the following decomposition of $T_n$ into a disjoint
union of subsets:
\begin{displaymath}
T_n=\bigcup_{A\subset\{2,3,\dots,n\}} T_A,
\end{displaymath}
where $T_A=\{\beta\in T_n: \beta\text{ contains the chain
}[\dots,1,a_1,\dots,a_k]\text{ and }\{a_1,\dots,a_k\}=A\}$. Set
$\overline{A}=N\setminus A$. If we substitute the chain
$[b_1,\dots,b_m,1,a_1,\dots,a_k]$ with its initial subchain
$[b_1,\dots,b_m]$, then every $\beta\in T_A$ can be transformed into
the element $\overline{\beta}$ from the set $\tilde{T}_A$ of all those
nilpotent elements from $\IS_{\overline{A}}$, which are not defined in
the point $1$. Moreover, every such nilpotent will be obtained exactly
$|A|!$ times. Hence $|T_A|=|A|!\cdot |\tilde{T}_A|$.

On the other hand, the set $\overline{A}$ is $\alpha$-invariant for 
every $\alpha\in E_A$, moreover, the restriction
$\alpha|_{\overline{A}}$ is a nilpotent element from $\tilde{T}_A$. 
Since the restriction $\alpha|_{\overline{A}}$ does not depend on the 
permutational part of $\alpha$, we get $|E_A|=|A|!\cdot |\tilde{T}_A|$. 

Therefore $|T_A|=|E_A|$ for all $A\subset \{2,3,\dots,n\}$ and hence
$l(n,0)=|E(n,0)|=|T_n|$.

\begin{remark}\label{rem3}{\rm
The equality $|T_A|=|E_A|$ can also be proved purely combinatorially,
using a bijection, analogous to that, constructed in
Remark~\ref{r2-3}. 
}
\end{remark}

\eqref{ff2}.  The orbit of $1$ under the action of $\alpha$ has length
$1$ if and only if $1$ is a fixed point of $\alpha$. The elements from
$\IS_n$, for which $1$ is a fixed point, are identified with
$\IS_{\{2,3,\dots,n\}}\simeq\IS_{n-1}$ in a natural way.

\eqref{ff3}. If the orbit of $1$ under the action of $\alpha$ has
length $k>1$, the element $\alpha$ has one of the following three
types: 
\begin{enumerate}[(I)]
\item\label{fff1} $\alpha=(1,a_2,\dots,a_k)\dots$. We have
$(n-1)\dots(n-k)|\IS_{n-k}|=[n-1]_{k-1}|\IS_{n-k}|$ elements of this
type. 
\item\label{fff2} $\alpha=[1,a_2,\dots,a_k]\dots$. We again have
$[n-1]_{k-1}|\IS_{n-k}|$ elements of this type.
\item\label{fff3} $\alpha=[b_1,\dots,b_m,1,a_2,\dots,a_k]\dots$. 
\end{enumerate} 
With every $\alpha$ of type \eqref{fff3} we associate the pair
$(\beta,[b_1,\dots,b_m])$, where $\beta\in
\IS_{N\setminus\{1,a_2,\dots,a_k\}}\simeq\IS_{n-k}$  is obtained from 
$\alpha$ substituting the chain
$[b_1,\dots,b_m,1,a_2,\dots,a_k]$ with the chain $[b_1,\dots,b_m]$. It 
is obvious that this map is a bijection to the set 
\begin{displaymath}
\{(\beta,l):\beta \in\IS_{N\setminus\{1,a_2,\dots,a_k\}},l
\text{ is a chain from the
chain decomposition of }\beta\}.
\end{displaymath}
The elements $a_2,\dots,a_k$ can be chosen in $[n-1]_{k-1}$ different
ways, and the pair $(\beta,l)$ in $L_{n-k}$ different ways.
Hence the number of elements of type \eqref{fff3} equals 
$[n-1]_{k-1}\cdot L_{n-k}$. Adding up the last three numbers we
obtained, we complete the proof of the theorem.
\end{proof}

\begin{corollary}\label{cor305}
$|T_n|$ equals the total number of partial injections from the set 
$\{2,3,\dots,n\}$ to the set $\{1,2,\dots,n\}$ (or from 
$\{1,2,\dots,n\}$ to $\{2,3,\dots,n\}$).
\end{corollary}

\begin{proof}
Follows from Theorem~\ref{t4}\eqref{ff1} and natural identification of
$E(n,0)$ with  partial injections from $\{2,3,\dots,n\}$ to 
$\{1,2,\dots,n\}$. Inverses for the later partial injections are
exactly partial injections from $\{1,2,\dots,n\}$ to 
$\{2,3,\dots,n\}$.
\end{proof}

\section{Nilpotent elements}\label{s4}

Some aspects of combinatorial properties of nilpotent elements in
$\IS_n$ were studied in \cite{GH2,GK1,GK2,GP}, however, the main
objects in these papers were not the elements from $T_n$ but rather
certain nilpotent subsemigroups in $\IS_n$, that is some special
subsets of $T_n$. The problem of calculating the cardinalities 
of such subsemigroups lead to interesting combinatorial
schemes, involving such classical combinatorial objects as Bell
numbers, Catalan numbers, Stirling numbers of the 2nd kind and
others. An overview of the results in this direction can be found in
\cite{GM2}. 

In this section we will investigate the combinatorial properties of
the set $T_n$ itself. A striking phenomenon we discover is a kind
of duality between the cardinalities of certain combinatorial
sets, associated with $\IS_n$ and $T_n$. This duality
will also appear in the next section and in the present section
it will be visible in most statements. However, we do not
have any satisfactory explanation for its existence.
We start with the theorem, which is in some sense dual to 
Theorem~\ref{t4}. We denote by $L^{(n)}$ the total number of chains 
in the chain decompositions of elements in $T_n$, and by $l^{n,k}$ the
total number of orbits of $1$ of length $k$ for the elements in $T_n$.

\begin{theorem}\label{t5}
\begin{enumerate}[(i)]
\item\label{t5en1} $l^{n,0}=|\IS_{n-1}|$.
\item\label{t5en2} $|\{\alpha\in T_n:\text{ the chain decomposition of
}\alpha\text{ contains the chain }[1]\}|=|T_{n-1}|$.
\item\label{t5en3} $l^{n,k}=[n-1]_{k-1}\cdot (L^{(n-k)}+|T_{n-k}|)$
for $1<k\leq n$.
\end{enumerate}
\end{theorem}

\begin{proof}
To prove \eqref{t5en1} we note that, by definition, $l^{n,0}$ is the
cardinality of the set $B=\{\alpha\in T_n:1\not\in \dom(\alpha)\}$.
The chain decomposition of every element from the set 
$B$ has the form $\alpha=[a_1,\dots,a_k,1]\dots$, where $k\geq 0$. 
Let us order the elements in $\{a_1,\dots,a_k\}$ in the increasing 
order: $a_{i_1}<a_{i_2}<\dots<a_{i_k}$. Note that the set 
$A=N\setminus \{a_1,\dots,a_k,1\}$ is $\alpha$-invariant, and define
$\overline{\alpha}\in \IS_{\{2,3,\dots,n\}}$ in the following way:
$\overline{\alpha}|_A=\alpha|_A$,
$\overline{\alpha}|_{\{a_1,\dots,a_k\}}=
\left(\begin{array}{ccc}a_{i_1}&\dots&a_{i_k}\\a_1&\dots&a_k
\end{array}\right)$. The map $\alpha\mapsto \overline{\alpha}$ is
obviously a bijection from $B$ to $\IS_{\{2,3,\dots,n\}}$ and the
statement follows.

\eqref{t5en2} is obvious.

To prove \eqref{t5en3} we partition the elements of the set
\begin{displaymath}
\{\alpha\in T_n:\text{ the orbit of }1\text{ under the action of }
\alpha\text{ has length }k\}
\end{displaymath}
into two classes, with respect to whether $1$ is a starting point of
some chain in the chain decomposition of $\alpha$ or not. The chain
decomposition of every element $\alpha$ from the first class has the
form $\alpha=[1,a_1,\dots,a_{k-1}]\dots$, where $a_1,\dots,a_{k-1}$
can be chosen in $[n-1]_{k-1}$ different ways, and all the other
chains of $\alpha$ define some nilpotent element from 
$\IS_{N\setminus\{1,a_1,\dots,a_{k-1}\}}$. Hence the first class
contains $[n-1]_{k-1}\cdot |T_{n-k}|$ elements.

The chain decomposition of every element from the second class has the
form $\alpha=[b_1,\dots,b_m,1,a_1,\dots,a_{k-1}]\dots$, where
$m>0$. The elements $a_1,\dots,a_{k-1}$ again can be chosen in 
$[n-1]_{k-1}$ different ways. If we now fix $a_1,\dots,a_{k-1}$, we
can associate the corresponding elements $\alpha$ to the pair 
$(\beta,[b_1,\dots,b_m])$, where $\beta\in 
\IS_{N\setminus\{1,a_1,\dots,a_{k-1}\}}$ is obtained from $\alpha$
substituting the chain $[b_1,\dots,b_m,1,a_1,\dots,a_{k-1}]$ with the
chain $[b_1,\dots,b_m]$. This defines, for fixed $a_1,\dots,a_{k-1}$,
a bijection from the set of all corresponding $\alpha$ to the set 
\begin{displaymath}
\{(\beta,l):\beta\text{ is a nilpotent element from  }
\IS_{N\setminus\{1,a_1,\dots,a_{k-1}\}}, l\text{ is a chain of 
}\beta\}.
\end{displaymath}
Hence the second class contains $[n-1]_{k-1}\cdot L^{(n-k)}$ elements.
\end{proof}

\begin{remark}\label{rem5}{\rm
The first parts of Theorems~\ref{t4} and \ref{t5} are completely
dual to each other. The last parts of these theorems are almost
dual, however, one could not expect a perfect duality for this 
statement as there are no orbits of length $1$ for nilpotent 
elements.
}
\end{remark}

\begin{theorem}\label{t6}
\begin{enumerate}
\item\label{t6en1} $|T_n|=\frac{1}{n}L_n$.
\item\label{t6en2} $|\IS_n|=\frac{1}{n+1}L^{(n+1)}$.
\end{enumerate}
\end{theorem}

\begin{proof}
The element $\alpha\in\IS_n$ can have some fixed points only in the
case, when the permutational part of $\alpha$ is not trivial, that 
is if $\alpha$ is not nilpotent. For every $\alpha\not\in T_n$ let 
$A_{\alpha}=\dom(\alpha^n)$ and $\overline{A_\alpha}=N\setminus 
A_{\alpha}$. Consider the set
\begin{displaymath}
M_{\alpha}=\{\beta\in\IS_n:\dom(\beta^n)=A_{\alpha}\text{ and
}\alpha|_{\overline{A_{\alpha}}}=\beta|_{\overline{A_{\alpha}}}\}. 
\end{displaymath}
Since the permutational parts of the elements from $M_{\alpha}$
correspond to all permutations in $S_{A_{\alpha}}$, it follows that
the average number of fixed points for  elements in $M_{\alpha}$
equals $1$. Since $M_{\alpha_1}=M_{\alpha_2}$ or $M_{\alpha_1}\cap 
M_{\alpha_2}=\varnothing$ for arbitrary $M_{\alpha_1}$ and
$M_{\alpha_2}$, the sets $M_{\alpha}$ form a partition of
$\IS_n\setminus T_n$ into a disjoint union of subsets. Hence the
total number $P_n$ of the fixed points equals $|\IS_n\setminus
T_n|$. Theorem~\ref{t1} now implies $\frac{1}{n}L_n=|\IS_n|-
|\IS_n\setminus T_n|=|T_n|$. This proves \eqref{t6en1}.

To prove \eqref{t6en2} it is enough to show that the cardinalities of
the sets 
\begin{gather*}
A=\{(x,\alpha):x\in\{1,2,\dots,n+1\},\alpha\in\IS_{\{1,2,\dots,n+1\}
\setminus \{x\}}\}\quad\quad\text{ and }\\
B=\{(\beta,l):\beta\in T_{n+1},l\text{ is a chain of }\beta\}
\end{gather*}
are the same. For this we define the map $f:A\to B$ in the following
way. Let $(x,\alpha)\in A$ and 
$\left(\begin{array}{ccc}a_1&\dots&a_k\\a_{i_1}&\dots&a_{i_k}
\end{array}\right)$ be the permutational part of $\alpha$. Assume that 
$a_1<a_2<\dots<a_k$ and set $f((x,\alpha))=(\beta,l)\in B$, where
$l=[a_{i_1},\dots,a_{i_k},x]$ and $\beta$ is obtained from $\alpha$ 
substituting the permutational part with $l$.

We define the map $g:B\to A$, $g:(\beta,l)\mapsto (x,\alpha)$ in the 
following way: if $l=[a_1,\dots,a_k,a_{k+1}]$, we set $x=a_{k+1}$ and
$\alpha$ is obtained from $\beta$ substituting $l$ with the
permutational part $\left(\begin{array}{ccc}a_{i_1}&\dots&a_{i_k}\\
a_1&\dots&a_k\end{array}\right)$, where
$a_{i_1}<a_{i_2}<\dots<a_{i_k}$ are elements $a_1,\dots,a_k$, written
with respect to the natural increasing order.

It is easy to check that $f$ and $g$ are mutually inverse to each
other, which implies $|A|=|B|$ and completes the proof.
\end{proof}

\begin{theorem}\label{t7}
\begin{enumerate}
\item\label{t7en1} $|T_n|=|\IS_{n-1}|+L_{n-1}$.
\item\label{t7en2} $|\IS_n|=|T_n|+L^{(n)}$.
\end{enumerate}
\end{theorem}

\begin{proof}
We start with \eqref{t7en1}.
According to the first part of Theorem~\ref{t4}, we have $|T_n|=|B|$,
where $B=\{\alpha\in\IS_n:1\not\in\dom(\alpha)\}$. We partition $B$
into two disjoint subsets $B_1=\{\alpha\in B:1\not\in\ran(\alpha)\}$
and $B_2=\{\alpha\in B:1\in\ran(\alpha)\}$. The elements of $B_1$ are
identified with the elements of $\IS_{\{2,3,\dots,n\}}\simeq
\IS_{n-1}$ in a natural way. Hence $|B_1|=|\IS_{n-1}|$.

The elements from $B_2$ have chain decomposition of the form 
$\alpha=[b_1,\dots,b_k,1]\dots$, where $k>0$. Sending every such
$\alpha$ to the pair $(\beta,[b_1,\dots,b_k])$, where $\beta\in 
\IS_{\{2,3,\dots,n\}}$ is obtained from $\alpha$ substituting the 
chain $[b_1,\dots,b_k,1]$ with the chain $[b_1,\dots,b_k]$, we get a
bijection from $B_2$ to the set $\{(\beta,l):\beta\in 
\IS_{\{2,3,\dots,n\}}, l\text{ is a chain of }\beta\}$. Hence
$|B_2|=L_{n-1}$. 

Now we prove \eqref{t7en2}.
Using the first part of Theorem~\ref{t5}, we can substitute $\IS_n$ by
$B=\{\alpha\in T_{n+1}:n+1\not\in\dom(\alpha)\}$. The chain
decomposition of every $\beta\in B$ contains the chain of the form
$[a_1,\dots,a_k,n+1]$, where $0\leq k\leq n$. For $k=0$ the
corresponding elements are identified with $T_n$ in a natural way,
hence the number of such elements in $|T_n|$. If $k>0$, we map the
element $\beta$ to the pair $(\overline{\beta},[a_1,\dots,a_k])$,
where $\overline{\beta}\in \IS_n$ is obtained from $\beta$ by 
substituting the chain $[a_1,\dots,a_k,n+1]$ by the chain
$[a_1,\dots,a_k]$. It is easy to see that this map is a bijection to
the set $\{(\alpha,l):\alpha\in T_n, l\text{ is a chain of
}\alpha\}$. Hence the number of such pairs equals $L^{(n)}$, which
completes the proof.
\end{proof}

\begin{remark}\label{rem6}{\rm
The first part of Theorem~\ref{t7} implies that nilpotent elements
form a substantial part of $|\IS_n|$, in particular, the inequality
$|T_n|>|\IS_{n-1}|$ is very rough.
}
\end{remark}

The following statement provides a precise connection between the
cardinalities of $\IS_n$ and $T_n$:

\begin{proposition}\label{p7}
\begin{displaymath}
|\IS_n|=\sum_{k=0}^n [n]_{k}|T_{n-k}|
=\sum_{k=1}^n [n-1]_{k-1}(n+k)|T_{n-k}|.
\end{displaymath}
\end{proposition}

\begin{proof}
The first equality follows from the fact that for a fixed $k>0$ the 
number of elements in $\IS_n$, which have stable rank $k$ equals 
$\binom{n}{k}\cdot k!\cdot |T_{n-k}|=[n]_k\cdot |T_{n-k}|$. 

To prove the second equality one shows, analogously to the proof of 
Proposition~\ref{p6}, that the average number
of fixed points in elements of stable rank $k>0$ is $1$. Moreover, 
the total number of points in the domains of the permutational parts 
of these elements equals $k\cdot \binom{n}{k}\cdot k!\cdot
|T_{n-k}|$. Using Theorem~\ref{t1} and Theorem~\ref{t2} we now get
\begin{displaymath}
|\IS_n|=\sum_{k=1}^n\binom{n}{k}\cdot k!\cdot |T_{n-k}|+
\frac{1}{n}\sum_{k=1}^n k\binom{n}{k}\cdot k!\cdot |T_{n-k}|=
\sum_{k=1}^n [n-1]_{k-1}(n+k)|T_{n-k}|.
\end{displaymath}
\end{proof}

\begin{corollary}\label{c4}
\begin{displaymath}
|T_{n}|=\sum_{k=1}^{n}k[n-1]_{k-1}|T_{n-k}|.
\end{displaymath}
\end{corollary}

\begin{proof}
Follows from the right equality of Proposition~\ref{p7}.
\end{proof}

\section{Various asymptotics}\label{s5}

\begin{lemma}\label{l1}
For every $n>1$ the following holds
\begin{enumerate}[(1)]
\item\label{enen1} $2n-1\geq |T_n|/|T_{n-1}|\geq n+1$, moreover, both
inequalities are strict for $n>2$,
\item\label{enen2} $2n> |\IS_n|/|\IS_{n-1}|> n+1$.
\end{enumerate}
\end{lemma}

\begin{proof}
To prove \eqref{enen1} we consider a chain $[a_1,\dots,a_k]$ from the
chain decomposition of some $\alpha\in T_{n-1}$. Inserting the point
$n$ on different places into this chain we obtain $k+1$ different
chains $[n,a_1,\dots,a_k]$, $[a_1,n,a_2,\dots,a_k]$, \dots, 
$[a_1,\dots,a_k,n]$. If we now perform this for every chain from the
chain decomposition of $\alpha$, we will get $(n-1)+\deff(\alpha)$
different nilpotent elements in $T_n$. One more nilpotent element is
obtained by adding the chain $[n]$ to $\alpha$. Since $1\leq
\deff(\alpha)\leq n-1$, we get 
\begin{equation}\label{eq4-1}
2n-1\geq (n-1)+\deff(\alpha)+1\geq n+1.
\end{equation}
Therefore for every $\alpha\in T_{n-1}$ we get at least $n+1$ and at
most $2n-1$ different elements from $T_n$. Certainly, performing this
construction for all elements from $T_{n-1}$ we will obtain all elements
from $T_{n}$, moreover, each element will be obtained only once. Hence
\begin{equation}\label{eq4-2}
(2n-1)\cdot |T_{n-1}|\geq  |T_{n}| \geq (n+1)\cdot|T_{n-1}|.
\end{equation}
If $n>2$, then the left inequality in \eqref{eq4-1} is strict for 
all $\alpha\in T_{n-1}$ such that $\deff(\alpha)=1$, and the right 
inequality in \eqref{eq4-1} is strict for all $\alpha\in T_{n-1}$ 
such that $\deff(\alpha)=n-1$. Hence both inequalities in
\eqref{eq4-2} are strict in this case.

The proof of \eqref{enen2} is analogous with the following
differences: one can insert the point $n$ in a cycle of length $k$ in 
$k$ different ways, one can add both the cycle $(a)$ and the chain
$[a]$ to the chain decomposition of $\alpha\in\IS_{n-1}$.
\end{proof}

If we recall that $T_n$ contains exactly $L'(n,k)$ and $\IS_n$
contains exactly $R_{n,n-k}$ elements of defect $k$, the proof of
Lemma~\ref{l1} immediately implies

\begin{lemma}\label{l2}
\begin{enumerate}
\item\label{l2it1} $\displaystyle
|T_{n+1}|=\sum_{k=1}^n (n+k+1)L'(n,k)$,
\item\label{l2it2} $\displaystyle
|\IS_{n+1}|=\sum_{k=0}^n (n+k+2)R_{n,n-k}=\sum_{k=0}^n (2n-k+2)R_{n,k}$.
\end{enumerate}
\end{lemma}

\begin{lemma}\label{l3}
If $1\leq k<\sqrt{n+1}-1$ then $L'(n,k+1)>L'(n,k)$, and if
$\sqrt{n+1}-1<k<n$ then $L'(n,k+1)<L'(n,k)$.
\end{lemma}

\begin{proof}
We have $\frac{L'(n,k+1)}{L'(n,k)}=\frac{n-k}{k(k+1)}$ and we have
\begin{displaymath}
\frac{n-k}{k(k+1)}>1\Leftrightarrow
k^2+2k-n>0\Leftrightarrow1\leq k<\sqrt{n+1}-1.
\end{displaymath}
\end{proof}

Using analogous arguments we can see that

\begin{lemma}\label{l4}
If $1\leq k<n+\frac{1}{2}-\sqrt{n+5/4}$ then $R_{n,k+1}>R_{n,k}$, 
and if $n+\frac{1}{2}-\sqrt{n+5/4}<k<n$ then $R_{n,k+1}<R_{n,k}$.
\end{lemma}

\begin{lemma}\label{l5}
\begin{displaymath}
\lim_{n\to\infty}\frac{n\cdot L'(n,3[\sqrt{n}])}
{\sqrt{n}\cdot L'(n,2[\sqrt{n}])}=0\quad\quad\text{ and }\quad\quad
\lim_{n\to\infty}\frac{n\cdot R_{n,n-3[\sqrt{n}]}}
{\sqrt{n}\cdot R_{n,n-2[\sqrt{n}]}}=0.
\end{displaymath}
\end{lemma}

\begin{proof}
To prove the first formula we set $m=[\sqrt{n}]$. Using the Stirling
formula for $n!$ we get
\begin{multline*}
\frac{n\cdot L'(n,3m)}{\sqrt{n}\cdot L'(n,2m)}=
\frac{n\cdot \frac{(n-1)!}{(n-3m)!(3m-1)!}\cdot \frac{n!}{(3m)!}}
{\sqrt{n}\cdot \frac{(n-1)!}{(n-2m)!(2m-1)!}\cdot \frac{n!}{(2m)!}}=\\
=\sqrt{n}\cdot\frac{\sqrt{2\cdot(2m-1)(n-2m)}}{\sqrt{3\cdot(3m-1)(n-3m)}} 
\cdot \exp(m)\cdot \left(\frac{n-2m}{n-3m}\right)^{n-3m}\cdot 
\left(\frac{2m-1}{3m-1}\right)^{2m-1}\times\\\times 
\left(\frac{n-2m}{3m-1}\right)^{m}\cdot \left(\frac{4}{27}\right)^{m}
\cdot\frac{1}{m^m}\cdot(1+o(1)).
\end{multline*}
But we have
\begin{gather*}
\frac{\sqrt{2\cdot(2m-1)(n-2m)}}{\sqrt{3\cdot(3m-1)(n-3m)}}=
\frac{2}{3}(1+o(1)),\quad\quad
\left(\frac{n-2m}{n-3m}\right)^{n-3m}=\exp(m)(1+o(1)),\\
\left(\frac{2m-1}{3m-1}\right)^{2m-1}=\left(\frac{2}{3}\right)^{2m-1}
\cdot\exp(-1/3)\cdot(1+o(1)),\\ 
\left(\frac{n-2m}{3m-1}\right)^{m}\leq 
\left(\frac{m^2}{3m-1}\right)^{m}=\left(\frac{m}{3}\right)^{m}
\cdot\exp(1/3)\cdot(1+o(1)).
\end{gather*}
Hence
\begin{displaymath}
\frac{n\cdot L'(n,3m)}{\sqrt{n}\cdot L'(n,2m)}\leq 
\sqrt{n}\cdot \left(\frac{16\exp(2)}{3^6}\right)^{m}\cdot(1+o(1)).
\end{displaymath}
As $\sqrt{n}=m(1+o(1))$ and $16\exp(2)/3^6<1$, we obtain
\begin{displaymath}
\lim_{n\to\infty}\frac{n\cdot L'(n,3[\sqrt{n}])}
{\sqrt{n}\cdot L'(n,2[\sqrt{n}])}=0.
\end{displaymath}

The proof of the second formula is analogous, using
$R_{n,k}=\binom{n}{k}^2\cdot k!$.
\end{proof}

\begin{theorem}\label{t8}
\begin{displaymath}
\lim_{n\to\infty}\frac{|T_{n+1}|}{(n+2)|T_n|}=
\lim_{n\to\infty}\frac{|\IS_{n+1}|}{(n+2)|\IS_n|}=1.
\end{displaymath}
\end{theorem}

\begin{proof}
From Lemma~\ref{l1} it follows that for all $n$ we have
\begin{displaymath}
\frac{|T_{n+1}|}{(n+2)|T_n|}\geq
\frac{(n+2)|T_{n}|}{(n+2)|T_n|}=1,\quad\text{and }\quad
\frac{|\IS_{n+1}|}{(n+2)|\IS_n|}>
\frac{(n+2)|\IS_{n}|}{(n+2)|\IS_n|}=1.
\end{displaymath}
Hence to prove the theorem it is enough to show that both sequences,
$\frac{|T_{n+1}|}{(n+2)|T_n|}$ and $\frac{|\IS_{n+1}|}{(n+2)|\IS_n|}$, 
are majorized by a sequence, which converges to $1$. For the sequence
$\frac{|T_{n+1}|}{(n+2)|T_n|}$ we have, using
Lemma~\ref{l2}\eqref{l2it1}, the following: 
\begin{multline}\label{t8eq1}
\frac{|T_{n+1}|}{(n+2)|T_n|}=
\frac{\displaystyle \sum_{k=1}^n(n+k+1)L'(n,k)}
{\displaystyle (n+2)\sum_{k=1}^n L'(n,k)}<\\<
\frac{\displaystyle \sum_{k<3[\sqrt{n}]}(n+3[\sqrt{n}]+1)L'(n,k)+
\sum_{k\geq 3[\sqrt{n}]}(2n+1)L'(n,k)}
{\displaystyle (n+2)\sum_{k=1}^n L'(n,k)}<\\<
\frac{n+3[\sqrt{n}]+1}{n+2}+\frac{2n+1}{n+2}\cdot
\frac{\displaystyle \sum_{k\geq 3[\sqrt{n}]}L'(n,k)}
{\displaystyle \sum_{k=1}^n L'(n,k)}.
\end{multline}
By Lemma~\ref{l3} we have
\begin{displaymath}
\sum_{k\geq 3[\sqrt{n}]}L'(n,k)<
\sum_{k\geq 3[\sqrt{n}]}L'(n, 3[\sqrt{n}])<n\cdot L'(n, 3[\sqrt{n}])
\end{displaymath}
and 
\begin{displaymath}
\sum_{k=1}^nL'(n,k)>
\sum_{[\sqrt{n}]\leq k\leq 2[\sqrt{n}]}L'(n,k)>[\sqrt{n}]\cdot 
L'(n, 2[\sqrt{n}]).
\end{displaymath}
Applying the first part of Lemma~\ref{l5} we get that the second
summand of \eqref{t8eq1} converges to $0$. It is obvious that the
first summand converges to $1$, which completes the proof for the
sequence $\frac{|T_{n+1}|}{(n+2)|T_n|}$.

For the sequence $\frac{|\IS_{n+1}|}{(n+2)|\IS_n|}$ we have, 
using Lemma~\ref{l2}\eqref{l2it2}, the following:
\begin{multline}\label{t8eq2}
\frac{|\IS_{n+1}|}{(n+2)|\IS_n|}=
\frac{\displaystyle \sum_{k=0}^n(2n-k+2)R_{n,k}}
{\displaystyle (n+2)\sum_{k=0}^n R_{n,k}}<\\<
\frac{\displaystyle \sum_{k<n-3[\sqrt{n}]}(2n+2)R_{n,k}+
\sum_{k\geq n-3[\sqrt{n}]}(n+3[\sqrt{n}]+2)R_{n,k}}
{\displaystyle (n+2)\sum_{k=1}^n R_{n,k}}<\\<
\frac{2n+2}{n+2}\cdot 
\frac{\displaystyle \sum_{k< n-3[\sqrt{n}]}R_{n,k}}
{\displaystyle \sum_{k=1}^n R_{n,k}}+\frac{n+3[\sqrt{n}]+2}{n+2}.
\end{multline}
By Lemma~\ref{l4} we have
\begin{displaymath}
\sum_{k<n-3[\sqrt{n}]}R_{n,k}<
\sum_{k<n-3[\sqrt{n}]}^nR_{n,n-3[\sqrt{n}]}<n\cdot R_{n,n-3[\sqrt{n}]}
\end{displaymath}
and 
\begin{displaymath}
\sum_{k=1}^n R_{n,k}>
\sum_{n-3[\sqrt{n}]\leq k\leq n-2[\sqrt{n}]}R_{n,k}>
[\sqrt{n}]\cdot R_{n,n-2[\sqrt{n}]}.
\end{displaymath}
Applying the second part of Lemma~\ref{l5} we get that the first
summand of \eqref{t8eq1} converges to $0$. It is obvious that the
second summand converges to $1$, which completes the proof.
\end{proof}

\begin{theorem}\label{t9}
\begin{displaymath}
\lim_{n\to \infty}\frac{|T_n|}{|\IS_n|}=0.
\end{displaymath}
\end{theorem}

\begin{proof}
Using the first statement of Theorem~\ref{t7}, Proposition~\ref{p5}
and Lemma~\ref{l2}\eqref{l2it2} we have
\begin{multline*}
\frac{|T_n|}{|\IS_n|}=\frac{|\IS_{n-1}|+L_{n-1}}{|\IS_n|}=
\frac{|\IS_{n-1}|}{|\IS_n|}+\\+
\frac{\displaystyle \sum_{k=0}^{n-1}(n-k-1)R_{n-1,k}}{|\IS_n|}=
\frac{|\IS_{n-1}|}{|\IS_n|}+
\frac{\displaystyle \sum_{k=0}^{n-1}(n-k-1)R_{n-1,k}}
{\displaystyle \sum_{k=0}^{n-1}(2n-k)R_{n-1,k}}.
\end{multline*}
By Lemma~\ref{l1}\eqref{enen2} the first summand of the last sum
converges to $0$. Let us study the second summand in more detail:
\begin{multline*}
\frac{\displaystyle \sum_{k=0}^{n-1}(n-k-1)R_{n-1,k}}
{\displaystyle \sum_{k=0}^{n-1}(2n-k)R_{n-1,k}}<
\frac{\displaystyle \sum_{k<n-3[\sqrt{n-1}]}(n-k-1)R_{n-1,k}+
\sum_{k\geq n-3[\sqrt{n-1}]}(n-k-1)R_{n-1,k}}
{\displaystyle \sum_{k\geq n-3[\sqrt{n-1}]}(2n-k)R_{n-1,k}}<\\<
\frac{\displaystyle n\cdot \sum_{k<n-3[\sqrt{n-1}]}R_{n-1,k}}
{\displaystyle n\cdot \sum_{k\geq n-3[\sqrt{n-1}]}R_{n-1,k}}+
\frac{\displaystyle 3[\sqrt{n-1}]\cdot 
\sum_{k\geq n-3[\sqrt{n-1}]}R_{n-1,k}}
{\displaystyle n\cdot \sum_{k\geq n-3[\sqrt{n-1}]}R_{n-1,k}}.
\end{multline*}
As in the proof of Theorem~\ref{t8}, Lemma~\ref{l4} and the second
part of Lemma~\ref{l5} guarantee that the first summand in the last
sum converges to $0$. It is obvious that the second summand
$3[\sqrt{n-1}]/n$ converges to $0$ as well. This completes the proof.
\end{proof}

\begin{theorem}\label{t10}
Let $m\in\N$ be fixed. Then the distribution of the ranks of the
elements of $\IS_n$ modulo $m$ is asymptotically uniform, that is
for all $p\in\Z$ we have
\begin{displaymath}
\lim_{n\to\infty}
\frac{\displaystyle \sum_{k\in A_{n,p}} 
R_{n,k}}{|\IS_n|}=\frac{1}{m},
\end{displaymath}
where $A_{n,p}=\{x\in\Z:0\leq x\leq n,\, x\equiv p\mod m\}$.
\end{theorem}

\begin{proof}
Denote $F_p= \sum_{k\in A_{n,p}} R_{n,k}$ and let $k_0=\lceil
n+1/2-\sqrt{n+5/4}\rceil$. For $p\in\Z$ let $A_{n,p}^0$ denote the set 
of all $x\in A_{n,p}$ satisfying $x\leq k_0$, and 
$A_{n,p}^1=A_{n,p}\setminus A_{n,p}^0$. From Lemma~\ref{l4} it follows 
that $R_{n,k}$ is increasing for $0\leq k\leq k_0$ and decreasing for
$k_0\leq k\leq n$. For $k_0$ the value $R_{n,k_0}$ is the maximal one
(for a fixed $n$). Hence for all $p,q$ we have
\begin{displaymath}
|F_p-F_q|\leq \left|\sum_{k\in A_{n,p}^0}R_{n,k}-
\sum_{k\in A_{n,q}^0}R_{n,k}\right|+
\left|\sum_{k\in A_{n,p}^1}R_{n,k}-
\sum_{k\in A_{n,q}^1}R_{n,k}\right|<2R_{n,k_0}.
\end{displaymath}

\begin{lemma}\label{l6}
Let $n>10000$ and $|k-k_0|<\frac{\sqrt[4]{n}}{6}-1$. 
Then $\left|\frac{R_{n,k+1}}{R_{n,k}}-1\right|<\frac{1}{\sqrt[4]{n}}$.
\end{lemma}

\begin{proof}
For  $s=k-(n+1/2-\sqrt{n+5/4})$ we obviously have 
$|s|<\frac{\sqrt[4]{n}}{6}$. By direct calculation we get
\begin{displaymath}
\left|\frac{R_{n,k+1}}{R_{n,k}}-1\right|=
\left|\frac{n^2-2nk+k^2-k-1}{k+1}\right|=
\left|\frac{s^2-2s\sqrt{n+5/4}}{n+3/2-\sqrt{n+5/4}+s}\right|.
\end{displaymath}
Again by direct calculation it is easy to  show that for
$|s|\leq 1$ we have
\begin{displaymath}
\left|\frac{s^2-2s\sqrt{n+5/4}}{n+3/2-\sqrt{n+5/4}+s}\right|<
\left|\frac{4\sqrt{n+5/4}}{n-\sqrt{n+5/4}}\right|<\frac{6\sqrt{n}}{n}<
\frac{1}{\sqrt[4]{n}};
\end{displaymath}
and that for $|s|> 1$ we have
\begin{displaymath}
\left|\frac{s^2-2s\sqrt{n+5/4}}{n+3/2-\sqrt{n+5/4}+s}\right|<
\left|\frac{4s\sqrt{n+5/4}}{n-2\sqrt{n}}\right|<\frac{6|s|}{\sqrt{n}}<
\frac{1}{\sqrt[4]{n}}.
\end{displaymath}
\end{proof}

\begin{lemma}\label{l7}
For all $n$ big enough the inequality 
$k_0-\left[\frac{\sqrt[4]{n}}{6}\right]+1\leq k\leq k_0+
\left[\frac{\sqrt[4]{n}}{6}\right]-1$ implies the inequality
$\frac{R_{n,k_0}}{R_{n,k}}<2$.
\end{lemma}

\begin{proof}
From Lemma~\ref{l6} it follows that for all such $k$ we have
\begin{displaymath}
\frac{R_{n,k_0}}{R_{n,k}}<
\left(1+\frac{1}{\sqrt[4]{n}}\right)^{\sqrt[4]{n}/6}
=e^{1/6}(1+o(1)).
\end{displaymath}
The remark that  $e^{1/6}<2$ completes the proof.
\end{proof}

From Lemma~\ref{l7} it follows that for all $n$ big enough and for all
$p$ and $q$ we have
\begin{displaymath}
\frac{|F_p-F_q|}{|\IS_n|}<\frac{2R_{n,k_0}}{|\IS_n|}<
\frac{2R_{n,k_0}}{\displaystyle \sum_{k\in B_n}R_{n,k}}
<\frac{2R_{n,k_0}}{2\left(\left[\frac{\sqrt[4]{n}}{6}\right]-1\right)
\cdot\frac{R_{n,k_0}}{2}}=
\frac{2}{\left[\frac{\sqrt[4]{n}}{6}\right]-1},
\end{displaymath}
where $B_n=\left\{k\in \Z: k_0-\left[\frac{\sqrt[4]{n}}{6}\right]+1\leq 
k\leq k_0+\left[\frac{\sqrt[4]{n}}{6}\right]-1\right\}$. Hence 
\begin{displaymath}
\lim_{n\to\infty}\frac{F_p}{|\IS_n|}=
\lim_{n\to\infty}\frac{F_q}{|\IS_n|}.
\end{displaymath}
As $F_0+F_1+\dots+F_{m-1}=|\IS_{n}|$, we finally get 
$\displaystyle \lim_{n\to\infty}\frac{F_p}{|\IS_n|}=\frac{1}{m}$.
\end{proof}

\section{Random products}\label{s6}

We consider the products $x_1x_2\dots x_k$ of elements from $\IS_n$
of length $k$. We assume that the elements $x_1,x_2,\dots, x_k$ are
chosen randomly and independently, with the uniform distribution of
probabilities, that is the probability to choose the element
$\alpha\in\IS_n$ does not depend on $\alpha$ and equals
$\frac{1}{|\IS_n|}$.  

\begin{lemma}\label{l8}
Given $\alpha\in\IS_n$, the probability of the following random event 
``the random product $x_1x_2\dots x_k$ of elements from $\IS_n$
of length $k$ equals $\alpha$'' depends only on $\rank(\alpha)$.
\end{lemma}

\begin{proof}
Let $\rank(\alpha)=\rank(\beta)=m$ and 
$\alpha=\left(\begin{array}{ccc}a_1&\dots&a_m\\
b_1&\dots&b_m\end{array}\right)$, 
$\beta=\left(\begin{array}{ccc}c_1&\dots&c_m\\
d_1&\dots&d_m\end{array}\right)$. Let $\mu,\tau\in S_n$ be such that 
$\mu=\left(\begin{array}{cccc}c_1&\dots&c_m&\dots\\
a_1&\dots&a_m&\dots\end{array}\right)$, 
$\tau=\left(\begin{array}{cccc}b_1&\dots&b_m&\dots\\
d_1&\dots&d_m&\dots\end{array}\right)$. Then $\mu\alpha\tau=\beta$ and 
the map
\begin{gather*}
\{(x_1,\dots,x_k):x_1\dots x_k=\alpha\}\to
\{(y_1,\dots,y_k):y_1\dots y_k=\beta\} ,\\
(x_1,\dots,x_k)\mapsto (\mu x_1,x_2,\dots,x_{k-1},x_k\tau)
\end{gather*}
is injective. Hence $\Pr(\alpha=x_1\dots x_k)\leq 
\Pr(\beta=y_1\dots y_k)$. The opposite inequality follows by
switching $\alpha$ and $\beta$.
\end{proof}

Let $P_{k,n}^{(i)}$ denote the probability of the random event
``the random product $x_1x_2\dots x_k$ of elements from $\IS_n$
of length $k$ is equal to a fixed element of rank $i$''. From
Lemma~\ref{l8} it follows that $P_{k,n}^{(i)}$ is well-defined, that
it does not depend on the choice of the element of rank $i$.

\begin{corollary}\label{cnew10101}
Let $M\subset\IS_n$ and $m_i$, $i=0,\dots n$, denote the number of
elements in $M$ of rank $i$. Then the probability of the following 
random event: ``the random product $x_1x_2\dots x_k$ of elements 
from $\IS_n$ of length $k$ belongs to $M$'' equals
$\displaystyle \sum_{i=0}^n m_i P_{k,n}^{(i)}$. In particular, 
the probability of the random event ``the random product 
$x_1x_2\dots x_k$ of elements from $\IS_n$ of length $k$ belongs to
$T_n$'' equals $\displaystyle \sum_{i=0}^{n-1} L'(n,n-i) 
P_{k,n}^{(i)}$.
\end{corollary}

\begin{proposition}\label{p8}
\begin{displaymath}
P_{k,n}^{(i)}=\left(\frac{|\IS_{n-i}|}{|\IS_n|}\right)^k\cdot
\left([n]_i\right)^{k-1}\cdot P_{k,n-i}^{(0)}.
\end{displaymath}
\end{proposition}

\begin{proof}
We fix $\alpha\in\IS_n$ such that $\rank(\alpha)=i$ and have
$P_{k,n}^{(i)}=\Pr(x_1\dots x_k=\alpha)$. 

Take any random product $x_1\dots x_k$ and set 
$A_0=\dom(x_1\dots x_k)$, $A_1=x_1(A_0)$, $A_2=x_2(A_1)$,\dots, 
$A_k=x_k(A_{k-1})=\ran(x_1\dots x_k)$. Set 
$B_j=N\setminus A_j$, $j=1,\dots,k$. Then with every $x_j$ we
associate two maps: the bijection $y_j=x_j|_{A_{j-1}}:A_{j-1}\to A_j$
and the partial injection $z_j=x_j|_{B_{j-1}}:B_{j-1}\to
B_j$. Moreover, the equality $x_1\dots x_k=\alpha$ becomes equivalent
to the following pair of equalities: $y_1\dots y_k=\alpha$, $z_1\dots
z_k=0$. The sets $A_1,A_2,\dots,A_{k-1}$ and bijections
$y_1,y_2\dots,y_{k-1}$ can be chosen arbitrarily and this can be done in
$\left([n]_i\right)^{k-1}$ different ways. After this choice the
factor $y_k$ is uniquely determined.

For every $j$, $j=0,1,\dots,k$, we fix a bijection
$B_j\to\{1,2,\dots,n-i\}$. Then every $z_j:B_{j-1}\to B_j$ is
associated in a natural way with a partial injection, 
$\hat{z}_j\in\IS_{n-i}$. Moreover, the condition $z_1\dots
z_k=0$ becomes equivalent to the condition $\hat{z}_1\dots
\hat{z}_k=0$. Since for the last equation the factors can be chosen in
$|\IS_{n-i}|^k\cdot P_{k,n-i}^{(0)}$ different ways, we get 
\begin{displaymath}
P_{k,n}^{(i)}=\left(\frac{|\IS_{n-i}|}{|\IS_n|}\right)^k\cdot
\left([n]_i\right)^{k-1}\cdot P_{k,n-i}^{(0)},
\end{displaymath}
which completes the proof.
\end{proof}

\begin{corollary}\label{c5}
\begin{displaymath}
\frac{|\IS_{n-i}|^k}{|\IS_n|^k}\cdot \left([n]_i\right)^{k-1}\geq 
P_{k,n}^{(i)}\geq
\frac{|\IS_{n-i}|^{k-1}}{|\IS_n|^k}\cdot \left([n]_i\right)^{k-1}.
\end{displaymath}
\end{corollary}

\begin{proof}
This follows from Proposition~\ref{p8} and the obvious inequality
$1\geq P_{k,n-i}^{(0)}\geq \frac{1}{|\IS_{n-i}|}$.
\end{proof}

\begin{corollary}\label{c6}
Let $n$ and $i>0$ be fixed. Then $\displaystyle
\lim_{k\to\infty}P_{k,n}^{(i)}=0$.
\end{corollary}

\begin{proof}
Using Corollary~\ref{c5} and Lemma~\ref{l1}\eqref{enen2} we get
\begin{multline*}
P_{k,n}^{(i)}\leq \frac{1}{[n]_i}
\left(\frac{|\IS_{n-i}|}{|\IS_n|}\cdot [n]_i\right)^k=\\=
\frac{1}{[n]_i}
\left(\frac{(n-i+1)|\IS_{n-i}|}{|\IS_{n-i+1}|}\cdot\dots\cdot
\frac{(n-1)|\IS_{n-2}|}{|\IS_{n-1}|}\cdot 
\frac{n\cdot|\IS_{n-1}|}{|\IS_{n}|}\right)^k<\\<
\frac{1}{[n]_i}
\left(\frac{(n-i+1)}{(n-i+2)}\cdot\dots\cdot
\frac{n-1}{n}\cdot 
\frac{n}{n+1}\right)^k=
\frac{1}{[n]_i}\left(\frac{n-i+1}{n+1}\right)^k.
\end{multline*}
But $\displaystyle
\lim_{k\to\infty}\frac{1}{[n]_i}\left(\frac{n-i+1}{n+1}\right)^k=0$,
which completes the proof.
\end{proof}

\begin{corollary}\label{c7}
Let $n$ be fixed. Then $\displaystyle
\lim_{k\to\infty}P_{k,n}^{(0)}=1$.
\end{corollary}

\begin{proof}
Since $x_1\dots x_k\in\IS_n$ we get $\displaystyle
\sum_{i=0}^n P_{k,n}^{(i)}\cdot R_{n,i}=1$. Since $R_{n,0}=1$, we
obtain the equality $\displaystyle P_{k,n}^{(0)}=1-\sum_{i=1}^n 
P_{k,n}^{(i)}\cdot
R_{n,i}$. From Corollary~\ref{c6} it follows that 
$\sum_{i=1}^n P_{k,n}^{(i)}\cdot R_{n,i}\to 0$ if $k\to\infty$,
and hence $P_{k,n}^{(0)}\to 1$.
\end{proof}

\begin{remark}\label{r7}{\rm 
Corollary~\ref{c7} implies that the semigroup $\IS_n$ is ``almost
nilpotent'' in the sense that for all $k$ big enough almost all
products $x_1\dots x_k$ of elements from $\IS_n$ equal $0$.
}
\end{remark}

\begin{corollary}\label{c8}
$P_{k,n}^{(n)}=\frac{(n!)^{k-1}}{|\IS_n|^k}$.
\end{corollary}

\begin{proof}
Follows from Proposition~\ref{p8} and the fact that $P_{k,0}^{(0)}=1$
as $\IS_0=\{0\}$.
\end{proof}

\begin{corollary}\label{c9}
For fixed $n$ and $i$ we have
\begin{displaymath}
P_{k,n}^{(i)}=\left(\frac{|\IS_{n-i}|}{|\IS_n|}\right)^k\cdot
\left([n]_i\right)^{k-1}\cdot (1+o(1)).
\end{displaymath}
\end{corollary}

\begin{proof}
Follows from Proposition~\ref{p8} and Corollary~\ref{c7}.
\end{proof}

For $i,j\in\N$ we denote by $I(i,j)$ the number of partial injections
from $\{1,2,\dots,i\}$ to $\{1,2,\dots,j\}$. It is obvious that
$I(i,i)=|\IS_i|$, $I(i,j)=I(j,i)$, and
\begin{displaymath}
I(i,j)=\sum_{k=0}^{\min(i,j)}\binom{i}{k}\binom{j}{k}k!.
\end{displaymath}

Consider the $(n+1)\times (n+1)$-matrix $\mathcal{A}=(A_{i,j})$,
$i,j=0,1,2,\dots,n$, where 
\begin{displaymath}
A_{i,j}=\begin{cases}
\binom{n-i}{j-i}\binom{n}{j}\cdot j!\cdot I(n-i,n-j),& 
\text{ if }i\leq j,\\
0,& \text{otherwise.}
\end{cases}
\end{displaymath}

\begin{theorem}\label{t11}
For all positive integers $n$ and $k$ we have the following equality
of vectors
\begin{displaymath}
\left(P_{k,n}^{(0)},P_{k,n}^{(1)},\dots,P_{k,n}^{(n)}\right)^{t}=
\frac{1}{|\IS_n|^k}\mathcal{A}^{k-1}\cdot (1,1,\dots,1)^{t}.
\end{displaymath}
\end{theorem}

\begin{proof}
We use induction in $k$ and note that the statement is obvious for
$k=1$. Let us now calculate $P_{k+1,n}^{(i)}=\Pr(x_1\dots
x_kx_{k+1})=\alpha$, where $\alpha\in\IS_n$ is a fixed element of rank
$i$. It is obvious that
\begin{displaymath}
P_{k+1,n}^{(i)}=\sum_{j=i}^n
\Pr(x_1\dots x_kx_{k+1}=\alpha\text{ and }\rank(x_1\dots x_k)=j).
\end{displaymath}
The product $x_1\dots x_k$ can be arbitrary, satisfying 
$\dom(x_1\dots x_k)\supset \dom(\alpha)$. Under the additional
assumption $\rank(x_1\dots x_k)=j$, we get that the product 
$x_1\dots x_k$ can have exactly $\binom{n-i}{j-i}\binom{n}{j}j!$
different values, where $\binom{n-i}{j-i}$ is the total number of
extensions of $\dom(\alpha)$ up to $\dom(x_1\dots x_k)$,
$\binom{n}{j}$ is the number of ways to choose $\ran(x_1\dots x_k)$
and $j!$ is the number of ways to construct a bijection from
$\dom(x_1\dots x_k)$ to $\ran(x_1\dots x_k)$.

For a fixed $x_1\dots x_k$ the action of $x_{k+1}$ on $\ran(x_1\dots
x_k)$ is uniquely defined, and the action of $x_{k+1}$ on 
$N\setminus \ran(x_1\dots x_k)$ can be arbitrary with the only
restriction $x_{k+1}(N\setminus \ran(x_1\dots x_k))\subset
N\setminus\ran(\alpha)$. Hence for fixed $x_1\dots x_k$ we have
exactly $I(n-j,n-i)$ possibilities to choose $x_{k+1}$. This implies
that 
\begin{displaymath}
\Pr(x_1\dots x_kx_{k+1}=\alpha\text{ and }\rank(x_1\dots x_k)=j)=
P_{k,n}^{(j)}\cdot \binom{n-i}{j-i}\binom{n}{j}j!\cdot
\frac{I(n-j,n-i)}{|\IS_n|},
\end{displaymath}
where $P_{k,n}^{(j)}\cdot \binom{n-i}{j-i}\binom{n}{j}j!$ is the
probability of the occurrence of the necessary factor $x_1\dots x_k$, 
and $\frac{I(n-j,n-i)}{|\IS_n|}$ is the probability of the occurrence
of the independent necessary factor $x_{k+1}$.

Therefore
\begin{displaymath}
P_{k+1,n}^{(i)}=\sum_{j=i}^n
P_{k,n}^{(j)}\cdot \binom{n-i}{j-i}\binom{n}{j}j!\cdot
\frac{I(n-j,n-i)}{|\IS_n|}=\frac{1}{|\IS_n|}\sum_{j=i}^n
P_{k,n}^{(i)}A_{i,j},
\end{displaymath}
and hence 
\begin{displaymath}
\left(P_{k+1,n}^{(0)},\dots,P_{k+1,n}^{(n)}\right)^{t}=
\frac{1}{|\IS_n|}\cdot \mathcal{A}\cdot 
\left(P_{k,n}^{(0)},\dots,P_{k,n}^{(n)}\right)^{t}.
\end{displaymath}
Taking into account the inductive assumption we complete the proof.
\end{proof}

We remark that the matrix $\mathcal{A}$ is upper triangular with the 
positive integers $A_{i,i}=[n]_i|\IS_{n-i}|$ on the diagonal. Hence these 
numbers are the eigenvalues of $\mathcal{A}$. Furthermore, according
to Lemma~\ref{l1}\eqref{enen2}, we have
\begin{displaymath}
\frac{A_{i,i}}{A_{i+1,i+1}}=
\frac{[n]_i|\IS_{n-i}|}{[n]_{i+1}|\IS_{n-i-1}|}=
\frac{|\IS_{n-i}|}{(n-i)|\IS_{n-i-1}|}>
\frac{n-i+1}{n-i}>1,
\end{displaymath}
and thus all eigenvalues of $\mathcal{A}$ are different. Hence 
$\mathcal{A}$ has $n+1$ linearly independent eigenvectors.

\begin{proposition}\label{p9A}
The vectors
\begin{displaymath}
\begin{array}{l}
f_0=(R_{n,0},0,\dots,0)^t,\\
f_1=(-R_{n,1},R_{n-1,0},0,\dots,0)^t,\\
\dots\\
f_k=((-1)^kR_{n,k},(-1)^{k-1}R_{n-1,k-1},\dots,R_{n-k,0},0,\dots,0)^t,\\
\dots\\
f_n=((-1)^nR_{n,n},(-1)^{n-1}R_{n-1,n-1},\dots,
-R_{1,1},R_{0,0})^t\\ 
\end{array}
\end{displaymath}
are the eigenvectors of $\mathcal{A}$ with eigenvalues $A_{0,0}$,
$A_{1,1}$,\dots, $A_{n,n}$ respectively.
\end{proposition}

\begin{proof}
We are going to prove the statement using induction in $n$. For this
we have to denote the matrix $\mathcal{A}$ of order $n+1$ by 
$\mathcal{A}_n$ and the corresponding vectors $f_0,\dots,f_n$ by
$f_0^{(n)},\dots,f_n^{(n)}$ respectively. Under this notation we have
\begin{displaymath}
\mathcal{A}_n=\left(
\begin{array}{c|ccc}
I(n,n) & R_{n,1}I(n,n-1) & \dots & R_{n,n}I(n,0)\\
\hline
 & & &  \\
0 &  & n\cdot\mathcal{A}_{n-1}  & \\
 & & &  
\end{array}
\right)
\end{displaymath}
and $f_{k}^{(n)}=((-1)^k R_{n,k}| f_{k-1}^{(n-1)})^t$.

For $n=0$ we have $\mathcal{A}_0=(1)$ and $f_0^{(0)}=(1)$ and the
statement is obvious.

Let us now assume that the statement is true for
$\mathcal{A}_{n-1}$. Then
\begin{multline*}
\mathcal{A}_n\cdot f_k^{(n)}=\left(
\begin{array}{c|ccc}
I(n,n) & R_{n,1}I(n,n-1) & \dots & R_{n,n}I(n,0)\\
\hline
 & & &  \\
0 &  & n\cdot\mathcal{A}_{n-1}  & \\
 & & &  
\end{array}
\right)\cdot
\left(\begin{array}{c}(-1)^k R_{n,k}\\ \hline \\ f_{k-1}^{(n-1)}\\ \\ 
\end{array}
\right)=\\=
\left(\begin{array}{c}
I(n,n)\cdot (-1)^k R_{n,k}+\left(
R_{n,1}I(n,n-1),\dots,R_{n,n}I(n,0)
\right)\cdot f_{k-1}^{(n-1)}
\\ \hline \\
n\mathcal{A}_{n-1}\cdot f_{k-1}^{(n-1)}\\ \\  
\end{array}
\right).
\end{multline*}

From the inductive assumption we get 
$\mathcal{A}_{n-1}\cdot f_{k-1}^{(n-1)}=[n-1]_{k-1}|\IS_{n-k}|f_{k-1}^{(n-1)}$
and hence $n\mathcal{A}_{n-1}\cdot
f_{k-1}^{(n-1)}=[n]_{k}|\IS_{n-k}|f_{k-1}^{(n-1)}$.

The only thing, which is left to complete the proof, is to show the
following equality for the first coordinate:
\begin{equation}\label{eq1pA}
\sum_{i=0}^k R_{n,i}I(n,n-i)\cdot (-1)^{k-i}\cdot
R_{n-i,k-i}=[n]_k I(n-k,n-k)\cdot (-1)^k R_{n,k}.
\end{equation}
But we have $R_{n,i}\cdot R_{n-i,k-i}=R_{n,k}\binom{k}{i}$, and hence,
canceling $R_{n,k}\cdot (-1)^{k}$, we reduce \eqref{eq1pA} to the
following equality:
\begin{equation}\label{eq2pA}
\sum_{i=0}^k (-1)^{i}\binom{k}{i} I(n,n-i)=[n]_k I(n-k,n-k).
\end{equation}

To prove \eqref{eq2pA} we count the number $F$ of those
$\alpha\in\IS_n$, for which $\dom(\alpha)\supset \{1,2,\dots,k\}$, in
two different ways. The number of those $\alpha\in\IS_n$, which
are not defined in $a_1,\dots,a_i$, equals $I(n-i,n)$. Therefore,
using the principle of inclusion and exclusion, we get 
\begin{displaymath}
F=\sum_{i=0}^k (-1)^i \binom{k}{i}I(n,n-i).
\end{displaymath}
On the other hand, if $\alpha\in\IS_n$ satisfies
$\{1,2,\dots,k\}\subset \dom(\alpha)$, we can choose the values for 
$\alpha$ on the elements from $\{1,2,\dots,k\}$ in 
$\binom{n}{k}\cdot k!=[n]_k$ different ways. If the action of 
$\alpha$ on $\{1,2,\dots,k\}$ is already defined, the extension to $N$
is naturally identified with a partial injection on 
$N\setminus\{1,2,\dots,k\}$, and thus can be performed in 
$I(n-k,n-k)$ different ways. Hence $F=[n]_k\cdot I(n-k,n-k)$, which
completes the proof of \eqref{eq2pA} and of the proposition.
\end{proof}

\begin{proposition}\label{p9B}
\begin{displaymath}
\sum_{k=0}^n (-1)^k|\IS_{n-k}|\cdot R_{n,k}=1.
\end{displaymath}
\end{proposition}

\begin{proof}
As we have seen in the proof of Proposition~\ref{p9A}, the number of
those $\alpha\in \IS_n$, which are defined in the given $k$ points,
equals $[n]_k\cdot I(n-k,n-k)=[n]_k\cdot |\IS_{n-k}|$. Hence, by the
principle of inclusion and exclusion, the number of those elements in 
$\IS_n$, which are not defined in any point, equals
\begin{displaymath}
\sum_{k=0}^n (-1)^k\binom{n}{k} [n]_k|\IS_{n-k}|=
\sum_{k=0}^n (-1)^k R_{n,k}|\IS_{n-k}|.
\end{displaymath}
On the other hand, $\IS_n$ contains exactly one element, $0$, which is
not defined in any point.
\end{proof}

\begin{corollary}\label{c9C}
The vector $(1,1,\dots,1)^t$ has coordinates
$(|\IS_n|,|\IS_{n-1}|,\dots,|\IS_1|,|\IS_0|)$ in the basis, formed 
by vectors $f_0,f_1,\dots,f_n$ (see Proposition~\ref{p9A}). 
\end{corollary}

\begin{proof}
The vectors $f_0,f_1,\dots,f_n$ form a basis as eigenvectors, which
correspond to different eigenvalues for a linear operator with simple
spectrum. Let $T=(t_{i,j})$ be the transformation matrix to the basis 
$f_0,f_1,\dots,f_n$. For the entries of this matrix we have:
\begin{displaymath}
t_{i,j}=\begin{cases}
(-1)^{j-i} R_{n-i,j-i}, & \text{ if } i\leq j,\\
0, & \text{ otherwise}.
\end{cases}
\end{displaymath}
The necessary statement is now equivalent to the equality
\begin{displaymath}
T\cdot (|\IS_n|,|\IS_{n-1}|,\dots,|\IS_1|,|\IS_0|)^t=(1,1,\dots,1)^t,
\end{displaymath}
which follows immediately from Proposition~\ref{p9B}.
\end{proof}
\vspace{1cm}

\begin{center}
\bf Acknowledgments
\end{center}

This paper was written during the visit of the first author to 
Uppsala University, which was supported by The Swedish Institute.
The financial support of The Swedish Institute and the hospitality 
of Uppsala University are gratefully acknowledged. For the second 
author the research was partially supported by The Swedish 
Research Council. We would like to thank the referee for useful
suggestions which led to the improvements in the paper.

\vspace{1cm}

\noindent
O.G.: Department of Mechanics and Mathematics, Kyiv Taras Shevchenko
University, 64, Volodymyrska st., 01033, Kyiv, UKRAINE,
e-mail: {\tt ganiyshk\symbol{64}univ.kiev.ua}
\vspace{0.5cm}

\noindent
V.M.: Department of Mathematics, Uppsala University, Box 480,
SE 751 06, Uppsala, SWEDEN, e-mail: {\tt mazor\symbol{64}math.uu.se},
web: ``http://www.math.uu.se/$\tilde{\hspace{2mm}}$mazor''
\vspace{0.5cm}


\begin{thebibliography}{99}
\bibitem[BH]{BH} B.Brown, P.M.Higgins, Finite full transformation 
semigroups as collections of random functions. Glasgow Math. J. 30 
(1988), no. 2, 203--211.
\bibitem[GK1]{GK1}
{O.G.Ganyushkin, T.V.Kormysheva}, On nilpotent subsemigroups of a 
finite symmetric inverse semigroup. (Russian) Mat. Zametki 56 (1994),
no. 3,29--35. Translation in Math. Notes 56 (1994), no. 3-4, 896--899 
(1995). 
\bibitem[GK2]{GK2}
{O.G.Ganyushkin, T.V.Kormysheva}, The structure of nilpotent 
subsemigroups of a finite inverse symmetric semigroup. (Ukrainian), Dopov.
Nats. Akad. Nauk Ukrainy 1995, no. 1, 8--10. 
\bibitem[GM1]{GM} O.Ganyushkin, V.Mazorchuk, The full finite Inverse 
symmetric semigroup $\IS_n$, Preprint 2001:37, Chalmers University of 
Technology and G{\"o}teborg University, G{\"o}teborg, 2001.
\bibitem[GM2]{GM2}
O.Ganyushkin, V.Mazorchuk, Combinatorics of nilpotents in  
$\IS_n$, Preprint 2002:11, Uppsala University, Sweden, 2002.
\bibitem[GP]{GP}
O.Ganyushkin, M.Pavlov, On the cardinalities of a class of nilpotent
semigroups and their automorphism groups, in ``Algebraic structures
and their applications'', Proceedings of the Ukrainian Mathematical 
Congress-2001, Kyiv, Institute of Mathematics NAS of Ukraine, 2002,
17-21.
\bibitem[GH1]{GH1} G.M.S.Gomes, J.M.Howie, On the ranks of certain 
semigroups of order-preserving transformations. Semigroup Forum 45 
(1992), no. 3, 272--282.
\bibitem[GH2]{GH2} G.M.S.Gomes, J.M.Howie, Nilpotents in finite 
symmetric inverse semigroups. Proc. Edinburgh Math. Soc. (2) 30
(1987), no. 3, 383--395.
\bibitem[Ha]{Ha} B.Harris, Probability distributions related to 
random mappings. Ann. Math. Statist. 31 (1960), 1045--1062. 
\bibitem[Hi1]{Hi1} P.M.Higgins, Combinatorial results for semigroups 
of order-preserving mappings. Math. Proc. Cambridge Philos. Soc. 113 
(1993), no. 2, 281--296.
\bibitem[Hi2]{Hi2} P.M.Higgins, Random products in semigroups of 
mappings. Lattices, semigroups, and universal algebra (Lisbon, 1988), 
89--99, Plenum, New York, 1990.
\bibitem[Hi3]{Hi3} P.M.Higgins, Techniques of semigroup theory. With a 
foreword by G. B. Preston. Oxford Science Publications. The Clarendon 
Press, Oxford University Press, New York, 1992.
\bibitem[Ho1]{Ho} J.M.Howie, Products of idempotents in certain 
semigroups of transformations. Proc. Edinburgh Math. Soc. (2) 17 
(1970/71), 223--236.
\bibitem[Ho2]{Ho1} J.M.Howie, An introduction to semigroup theory. 
L.M.S. Monographs, No. 7. Academic Press, London-New York, 1976.
\bibitem[Ka]{Ka} L.Katz, Probability of indecomposability of a random 
mapping function. Ann. Math. Statist. 26, (1955). 512--517.
\bibitem[Kr]{Kr} M.Kruskal, The expected number of components under a 
random mapping function. Amer. Math. Monthly 61, (1954), 392--397.
\bibitem[La]{La} M.V.Lawson, Inverse semigroups. The theory of partial 
symmetries. World Scientific Publishing Co., Inc., River Edge, NJ, 1998.
\bibitem[Li]{Li} S.Lipscomb, Symmetric inverse semigroups. Mathematical 
Surveys and Monographs, 46. American Mathematical Society, Providence, 
RI, 1996.
\bibitem[Pe]{Pe} M.Petrich, Inverse semigroups. Pure and Applied 
Mathematics. A Wiley-Interscience Publication. John Wiley \& Sons,
Inc., New York, 1984
\end{thebibliography}
\end{document}